\documentclass{article}
\usepackage[utf8]{inputenc}
\usepackage{graphicx}	\usepackage{amsmath,amssymb,amsfonts,amsthm}
\usepackage{xcolor}
\usepackage[pagewise]{lineno}
\begin{document}

		\title{\textbf{The Inverse first passage time method for a two dimensional Ornstein Uhlenbeck process with neuronal application}}

	\author{$\text{A. Civallero}$, $\text{C. Zucca}$\\
		Dept. of Mathematics ``G. Peano'', \\
		\footnotesize University of Torino, Turin, Italy\\
	}
\date{}
	\maketitle

	\begin{abstract}

The Inverse First Passage time problem seeks to determine the boundary corresponding to a given stochastic process and a fixed first passage time distribution. 
Here, we determine the numerical solution of this problem in the case of a two dimensional Gauss-Markov diffusion process. We investigate the boundary shape corresponding to Inverse Gaussian or Gamma first passage time distributions for different choices of the parameters, including heavy and light tails instances. Applications in neuroscience  framework are illustrated.

	\end{abstract}
	
	{{\it Keywords:} Inverse First-passage-time problem, two-dimensional Ornstein Uhlenbeck process, two-compartment leaky integrate and fire model, Gamma, Inverse Gaussian.}

\section{Introduction}
\label{intro}
In many situations arising from applications (i.e. neuroscience, finance, reliability, ...), the quantity of interest is the first time that a random quantity crosses a given fixed level. In a mathematical framework, this corresponds to the first passage time (FPTs) of a stochastic process through an eventually time dependent boundary. However, it can happen that the FPT distribution is known as well as the random process, while one is interested in determining the corresponding time dependent boundary. This is the so-called Inverse FPT problem.
This problem has been investigated both from a theoretical \cite{CCCS, EJ} and an empirical point of view \cite{Ab,LSZ,SVZ,ZS} in the one-dimensional case. In \cite{CCCS} the existence and uniqueness of the solution of the Inverse FPT is studied.  In \cite{EJ} the problem is interpreted in terms of an optimal stopping problem. A numerical algorithm has been proposed in \cite{ZS} for the Wiener process. In \cite{SVZ} the Inverse FPT of an Ornstein Uhlenbeck (OU) process has been studied and it has been applied to a classification method with applications to neuroscience. The same framework is in \cite{LSZ}, where possible thresholds corresponding to Gamma distributed FPTs for an OU process has been investigated with modelling purposes. 

Here, we resort again to the Inverse FPT method, we generalize the algorithms to a two-dimensional OU process and we study the possibility to have Inverse Gaussian (IG) or Gamma distributed FPTs. The choice of these two distributions grounds on their role in neurosciences \cite{GM,KYP,LR,RS,SG,TIF, Y} and reliability theory \cite{FC,IP}. 

In Section \ref{Sec: model} we introduce the two dimensional Gauss-Markov process of interest, underlying some properties that we will use to deal with the Inverse FPT algorithm. In Section \ref{Sec: IFPT} we introduce the Inverse FPT method for a two dimensional process but we postpone to a future work the mathematical discussion about its convergence. In Section \ref{Sec:2 ISI} we apply the algorithm to two choices of the FPT distribution, determining the thresholds corresponding to IG or Gamma FPTs probability density function (pdf). We underline the differences between the two models and we explain how heavy or light tails influence the boundary behavior. 
The last section discusses the obtained results in neuroscience contest. The two compartment model of Leaky Integrate and Fire type  presented in \cite{LR} and studied in \cite{BSZ} describes the membrane potential evolution of a neuron as a two-dimensional OU process. Hence, in Section \ref{Sec: Results} we reinterpret the Inverse FPT results in this framework.



\section{The two dimensional Ornstein Uhlenbeck process} \label{Sec: model}



Let us consider a stochastic process $X=\{(X_1(t),X_2(t)), t\geq 0\}$ that is solution of the following stochastic differential system
\begin{equation}\label{model}
\left\{\begin{array}{lll} dX_1(t)=\left\{-\alpha X_1(t) + \beta \left[X_2(t)-X_1(t)\right]\right\}dt\\
                            \vspace{0.1mm}\\
                            dX_2(t)=\left\{-\alpha X_2(t) + \beta \left[X_1(t)-X_2(t)\right] + \mu\right\}dt + \sigma dB_t \end{array} \right.
\end{equation}
with $X(0)=0$ and where $B$ is a one-dimensional standard Brownian motion. Here, $\alpha>0$, $\beta$, $\mu$ and $\sigma>0$ are constants. 

\noindent
To solve the stochastic differential system (\ref{model}), we rewrite it in matrix form  
\begin{eqnarray}\label{SDE}
dX(t)=[AX(t) + M(t)]dt + G dB(t),
\end{eqnarray}
 where 
$$A=\left(\begin{array}{cc} -\alpha-\beta & \beta \\ \beta & -\alpha-\beta \end{array}\right) \mbox{, } M(t)=M=\left(\begin{array}{ll} 0 \\ \mu \end{array}\right) \mbox{ and }  
G=\left(\begin{array}{cc} 0 & 0 \\ 0 & \sigma \end{array}\right) \mbox{.}$$
It is an autonomous linear stochastic differential equation, in particular it is a two-dimensional Ornstein Uhlenbeck  process, special case of a Gauss-Markov diffusion process \cite{A}.
The solution of (\ref{SDE}) is
\begin{eqnarray}\label{SDE solution}
\left\{ 
\begin{array}{l}
X_1(t)=\frac{\mu}{2}\left(\frac{1-e^{-\alpha t}}{\alpha}-\frac{1-e^{-(\alpha +2\beta)t}}{\alpha +2\beta}\right)+\frac{\sigma }{2}\int_0^t \left( e^{-\alpha(t-s)}-e^{-(\alpha+2\beta)(t-s)}\right)dB(s)\\
X_2(t)=\frac{\mu}{2}\left(\frac{1-e^{-\alpha t}}{\alpha}+\frac{1-e^{-(\alpha +2\beta)t}}{\alpha+2\beta}\right)+\frac{\sigma}{2} \int_0^t \left( e^{-\alpha(t-s)}+e^{-(\alpha+2\beta)(t-s)}\right)dB(s).
\end{array}
\right.
\end{eqnarray}
It is a Gaussian vector with mean
\begin{equation}\label{media}
m(t)=\mathbb{E}(X(t))=\left[ 
\begin{array}{l}
\frac{\mu}{2}\left(\frac{1-e^{-\alpha t}}{\alpha}-\frac{1-e^{-(\alpha +2\beta)t}}{\alpha +2\beta}\right)\\
\frac{\mu}{2}\left(\frac{1-e^{-\alpha t}}{\alpha}+\frac{1-e^{-(\alpha +2\beta)t}}{\alpha+2\beta}\right)
\end{array}
   \right]
\end{equation}
and variance-covariance matrix $Q(t-s)$ where,
 
\begin{equation}
Q(t)=\left[
  \begin{array}{ll}
  Q^{(11)}  &Q^{(12)}\\
  Q^{(12)}  &Q^{(22)}
  \end{array}
  \right](t)
 \end{equation}
and
\begin{eqnarray*}
Q^{(11)}(t)&=&\frac{1}{2}
 \left(\frac{1}{\alpha}-\frac{2}{\alpha+\beta}+\frac{1}{\alpha+2\beta}-e^{-2\alpha t}\left(\frac{1}{\alpha}-\frac{2e^{-2\beta t}}{a+\beta}+\frac{e^{-4\beta t}}{\alpha+2\beta} \right)\right)\\ 
 Q^{(12)}(t)&=&\frac{1}{2}\frac{1-e^{-2\alpha t}}{\alpha}-\frac{1-e^{-2(\alpha +2\beta)t}}{\alpha+2\beta}\\
 Q^{(22)}(t)&=&\frac{1}{2}\left(\frac{1}{\alpha}+\frac{2}{\alpha+\beta}+\frac{1}{\alpha+2\beta}-e^{-2\alpha t}\left(\frac{1}{\alpha}+\frac{2e^{-2\beta t}}{\alpha+\beta}+\frac{e^{-4\beta t}}{\alpha+2\beta} \right)  \right).
\end{eqnarray*}

\begin{figure}[t!]
\includegraphics[width=0.90\textwidth]{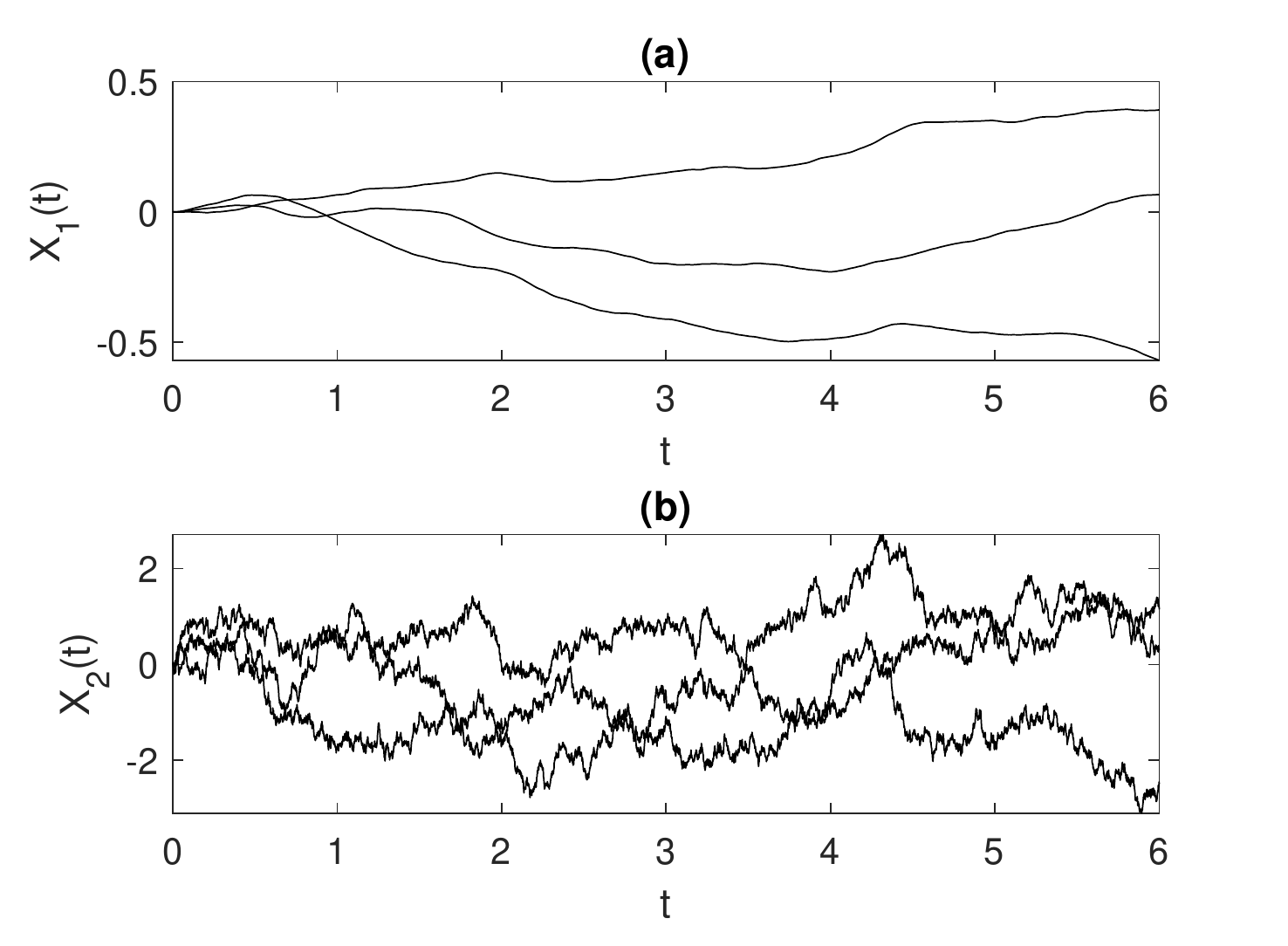}
	\caption{\label{traiettoria}Three sample paths of the two components of the process $X$. The parameters of the two compartment model are $\alpha=0.33$,
$\beta=0.2$, $\mu=0$ and $\sigma = 1$.}
\end{figure}

Trajectories of the process are plotted in Figure \ref{traiettoria}. The different behavior of the two components is evident: the noisy behavior is prevalent in $X_2$, while the first component $X_1$ is smoother. Indeed, as shown in (\ref{SDE solution}), the multiplying function of the random term on the first component reduces the noise effect.

We consider the first passage time of the first component of the process (\ref{model})
\begin{equation}
T=\inf\{t> 0: X_1(t)>S(t)\}
\end{equation}
where $S(t)$ is a continuous function with $S(0)\geq X_1(0)=0$. 

Note that it is possible to rewrite (\ref{SDE solution}) in iterative form. This version is useful for simulation purposes, in order to generate the trajectories in an exact way.  
Discretizing the time interval $[0,T]$ with the partition $\pi: 0=t_0 < t_1 < \dots <t_N =T$ in $N$ subintervals of constant length $h=\frac{T}{N}$, we can express the position of the process at time $t_{k+1}$ in terms of the position of the process at time $t_k$
  \begin{eqnarray}
  X(t_{k+1})&=&\frac{1}{2}\left[
  \begin{array}{ll}
  e^{-\alpha h}+e^{-(\alpha+2\beta)h} & e^{-\alpha h}-e^{-(\alpha+2\beta)h}\\
  e^{-\alpha h}-e^{-(\alpha+2\beta)h} & e^{-\alpha h}+e^{-(\alpha+2\beta)h}
\end{array}   \right]X(t_k)\\
&&+\frac{\mu}{2}\left[
  \begin{array}{l}
  \frac{1-e^{-\alpha h}}{\alpha}-\frac{1-e^{-(\alpha+2\beta)h}}{\alpha+2\beta}\\
  \frac{1-e^{-\alpha h}}{\alpha}+\frac{1-e^{-(\alpha+2\beta)h}}{\alpha+2\beta}\end{array}   \right] +\frac{\sigma}{2}I_k \nonumber
  \end{eqnarray}
where the term
  \begin{eqnarray}
  I_{k}&=&\left[
  \begin{array}{l}
\int_{t_{k}}^{t{k+1}}\left(  e^{-\alpha(t_{t+1}-s)}-e^{-(\alpha+2\beta)(t_{k+1}-s)}\right) dB(s)\\
\int_{t_{k}}^{t{k+1}}\left(  e^{-\alpha(t_{t+1}-s)}-e^{-(\alpha+2\beta)(t_{k+1}-s)}\right) dB(s)
\end{array}   \right],
  \end{eqnarray}
 known as innovation, is a Gaussian vector with zero mean and variance-covariance matrix $Q(h)$.


In the following we will also need the conditioned mean of the first component \begin{eqnarray}\label{media processo}
m^{(1)}(t|(X_1(\theta), X_2(\theta)),\theta) &=& \mathbb{E}(X_1(t)|X(\theta)=(X_1(\theta), X_2(\theta)))\\
&=&\frac{\mu}{2}
 \left(\frac{2\beta - \alpha e^{-\alpha (t-\theta)} -2\beta e^{-\alpha (t-\theta)} + \alpha e^{-(\alpha +2\beta ) (t-\theta)}}{\alpha(\alpha+2\beta)}\right)  \nonumber\\
&\quad &-\frac{X_1(\theta)e^{-\alpha(t-\theta)}}{2}(1+e^{-2\beta(t-\theta)}) - \frac{X_2(\theta)e^{-\alpha(t-\theta)}}{2} (1-e^{-2\beta(t-\theta)})\nonumber
\end{eqnarray}
and the conditioned variance of the first component 
\begin{eqnarray}
&&Q^{(11)}(t|(X_1(\theta),X_2(\theta)),\theta) = Var(X_1(t)|X(\theta)=(X_1(\theta), X_2(\theta)))\\
&&\quad \quad= \frac{\sigma^2 e^{-2\alpha (t-\theta)}}{8} \times \nonumber\\
&&\qquad \qquad\times \frac{2\alpha e^{-2\beta (t-\theta)}(\alpha +2\beta) -\alpha e^{-4\beta (t-\theta)}(\alpha+\beta) + 2\beta^2e^{2\alpha (t-\theta)} -\alpha^2 -3\alpha \beta -2\beta^2}{\alpha(\alpha +\beta)(\alpha +2\beta)}\nonumber
\end{eqnarray}


In some instances can be useful to transfer the time dependency from the boundary shape $S(t)$ to an input $M(t)$. Mathematically it is possible to relate these two situations with a simple space transformation.
Indeed, the space transformation
\begin{equation}\label{trasformazione}
Y_1(t) = X_1(t)-S(t)+\Sigma.
\end{equation}
changes our process $X$ given by (\ref{model}), originated in $X(0)=x_0$ in presence of a time dependent boundary $S(t)$
\begin{equation}
\left\{ 
\begin{array}{lll} 
dX_1(t)=\left\{-\alpha X_1(t) + \beta \left[X_2(t)-X_1(t)\right]\right\}dt\\
                            dX_2(t)=\left\{-\alpha X_2(t) + \beta \left[X_1(t)-X_2(t)\right] + \mu\right\}dt + \sigma dB_t\\
                     \vspace{0.1mm}\\
       
X(0)=x_0\\
S(t)                        
                             \end{array}
\right.
\end{equation}
into a two dimensional process characterized by time dependent
input $M(t)$ and constant threshold $\Sigma$
\begin{equation}
\left\{ 
\begin{array}{lll} 
dY_1(t)=\left\{-\alpha Y_1(t) + \beta \left[X_2(t)-Y_1(t)\right]+\mu_1(t)\right\}dt\\
dX_2(t)=\left\{-\alpha X_2(t) + \beta \left[Y_1(t)-X_2(t)\right] + \mu_2(t)\right\}dt + \sigma dB_t\\
                     \vspace{0.1mm}\\    
X(0)=x_0-S(0)+\Sigma\\
\Sigma.
     \end{array}
\right.
\end{equation}
Here, the term 
\begin{eqnarray}\label{M}
M(t)=\left[
\begin{array}{cc}
\mu_1(t)\\
\mu_2(t)
\end{array}
\right]=
\left[
\begin{array}{cc}
-(\alpha+\beta)(S(t)-\Sigma)-S'(t)\\
\mu+\beta(S(t)-\Sigma)
\end{array}
\right]
\end{eqnarray}
can be interpreted as an external input acting with different weights on the two compartments. 
Note that in (\ref{M}) $S(t)$ should be interpreted as a function of time and not as a boundary of a FPT problem. Indeed, in this case the boundary of the model is constant.

The stochastic process (\ref{model}) can be used to describe a system whose behavior depends by two components that are strictly correlated by the parameter $\beta$. The second component is driven by a random Gaussian noise and its evolution is stopped when the first component reaches a given fixed level. This model, known as two compartment model, has many interesting applications, for example in neuroscience, reliability and finance.

\section{Inverse first passage time method}
\label{Sec: IFPT}
The inverse FPT problem consists in searching the unknown boundary $S(t)$ given that the FPT density $f_T(t)$ is known.
We work under the assumption that the boundary $S(t)$ exists, it is unique and sufficiently regular. 

Let us consider a diffusion process $X=\{(X_1(t),X_2(t)), t\geq 0\}$, solution of the  stochastic differential equation (\ref{SDE}).
The proposed method is based on the numerical approximation of the following Volterra  integral equation \cite{BSZ}
\begin{equation} \label{equazione integrale_diffusione_inverso}
\begin{aligned}
1-Erf \Bigg( & \frac{S(t)-m^{(1)}(t)}{\sqrt{2Q^{(11)}(t)}} \Bigg) = \int_0^t d\theta f_T(\theta) \cdot \\
& \cdot \mathbb{E}_{Z(\theta)} \Bigg[ 1-Erf \Bigg( \frac{S(t)-m^{(1)}(t|(S(\theta),X_2(\theta)),\theta)}{\sqrt{2Q^{(11)}(t|(S(\theta),X_2(\theta)),\theta)}} \Bigg) \Bigg]
\end{aligned}
\end{equation}
where $Z(t)$ is a random variable that represents the position of the second component $X_2$ of the process when the first component $X_1$ hits the boundary at time $t$, i.e.
\begin{equation} \label{Z}
\mathbb{P}(Z(t)<z)=\mathbb{P}(X_2(T)<z|T=t,X(t_0) = y).
\end{equation} 

Let us fix a time interval $[0,\Theta]$ and a partition $\pi: 0=t_0 < t_1 < \dots <t_N =\Theta$ in $N$ subintervals of constant length $h=\frac{\Theta}{N}$. Using Euler formula for integrals \cite{AS}, equation (\ref{equazione integrale_diffusione_inverso}) can be approximated as \begin{equation} \label{appross_S}
\begin{aligned}
1- Erf \Bigg( & \frac{S^*(t_i)-m^{(1)}(t_i)}{\sqrt{2Q^{(11)}(t_i)}} \Bigg) = h \cdot \sum_{j=1}^i f_T(t_j) \cdot \\
& \cdot \mathbb{E}_{Z(t_j)} \Bigg[ 1-Erf \Bigg( \frac{S^*(t_i)-m^{(1)}(t_i|(S^*(t_j),X_2(t_j)),t_j)}{\sqrt{2Q^{(11)}(t_i|(S^*(t_j),X_2(t_j)),t_j)}} \Bigg) \Bigg]
\end{aligned}
\end{equation}
$\forall i = 1 \dots N$, \\
 
Equation \eqref{appross_S} represents a non linear system of  $N$ equations in $N$ unknown $S^*(t_1), \dots , S^*(t_N)$ that can be solved by means of root finding iterative algorithms \cite{At}. Its solution gives an approximation $S^*(t)$ of the boundary $S(t)$ in the partition points $\pi$. Note that in step $i$ the only unknown quantity is $S(t_i)$ and it is 
estimated using the boundary approximations $S^*(t_1), \dots , S^*(t_{i-1})$, computed in the previous steps.

The quantity
\begin{equation}\label{theta}
\theta_{i,k}=\mathbb{E}_{Z(t_k)} \Bigg[ 1-Erf \Bigg( \frac{S^*(t_i)-m^{(1)}(t_i|S^*(t_k|(S^*(t_k),X_2(t_k)),t_k)}{\sqrt{2Q^{(11)}(t_i|(S^*(t_k),X_2(t_k)),t_k)}} \Bigg) \Bigg]
\end{equation}
is not easily handled because it depends on the unknown time dependent boundary. In general, the 
computation of $\theta_{k,k} $ is not trivial 
but, performing a suitable limit on the considered process we can show that $\theta_{k,k}=2$ for each value of $k$. To compute (\ref{theta}) when  $k\neq i$, we use a Monte Carlo method: we simulate the process $ X $ until the first component exceeds the threshold and we save the corresponding value of $Z$. 
At step $i$, we need to compute $\theta_{i,k}$ for $k=1, \dots, i-1$. 
The presence of an expectation with respect to $Z(t_k)$ determines a difficulty for the estimation of (\ref{theta}) through Monte Carlo because we need the value of $Z(t_k)=X_2(t_k)$  at time $T=t_k$. To circumvent this problem we introduce an approximate approach as follows.
At step $i$ we approximate the threshold 
with a piecewise linear curve with knots in the already computed boundary values. Hence, for $\tau\in  [ t_ {j-1}, t_j]$, $\quad j=1,\dots, i-1$ we substitute the exact boundary with  
\begin{equation}\label{Barr appr}
\hat{S}(\tau)= \frac{S^*(t_j)-S^*(t_{j-1})}{t_j - t_{j-1}} \tau + \frac{t_j S^*(t_{j-1}) - S^*(t_j)t_{j-1}}{t_j - t_{j-1}}
\end{equation} 
and we simulate the process up to $ t_{i-1}$ or until it reaches the threshold. To compute $\theta_{i,k}$, $k=1,\cdots, i-1$ we use only the trajectories that crossed the approximated  boundary (\ref{Barr appr}) in a neighbourhood of $t_k$. Then, in correspondence to each of these sample paths we identify with $\{Z_k, k = 1 \dots M \}$  the sequence of values of the second component of the process $X$ (when the first component has exceeded the threshold). In this way, the Monte Carlo estimate for $ \theta_{i, j}$ is 
\begin{displaymath}
\tilde{\theta}_{i,j} = 1- \frac{\sum_{k=1}^M Erf \Bigg( \frac{S(t_i)-m^{(1)}(t_i|(S(t_j),Z_k),t_j)}{\sqrt{2Q^{(11)}(t_i|(S(t_j),Z_k),t_j)}} \Bigg)}{M}.
\end{displaymath}
It is possible to prove that this further approximation does not seriously influence the reliability of the algorithm.

\section{Examples}\label{Sec:2 ISI}
In this Section we illustrate the use of the Inverse FPT method through two examples.
The first situation
concerns FPTs with Inverse Gaussian distribution. The second one deals with the Gamma distribution. Lastly, a comparison between boundaries and drift terms arising in the two examples is developed. 

\subsection{Inverse Gaussian random variable}\label{Sec:IG}
The IG random variable $T$ has pdf
\begin{equation}\label{f_TIG}
f_T(t)=\left[ \frac{\lambda}{2\pi t^3}\right]^{1/2} \exp\left[ -\frac{\lambda(t-\rho)^2}{2\rho^2 t}\right], \quad \quad t\geq 0, 
\end{equation}
where $\rho>0$ is the mean and $\lambda>0$ is the shape parameter. 
Mean, variance and coefficient of variation are given by
\begin{eqnarray}\label{CV_IG}
\mathbb{E}(T)&=&\rho\\\nonumber
Var(T)&=&\frac{\rho^3}{\lambda}\\
CV&=&CV(T)=\frac{\sqrt{Var(T)}}{\mathbb{E}(T)}=\sqrt{\frac{\rho}{\lambda}}.\nonumber
\end{eqnarray}

\begin{figure}[t!]
\includegraphics[width=1\textwidth]{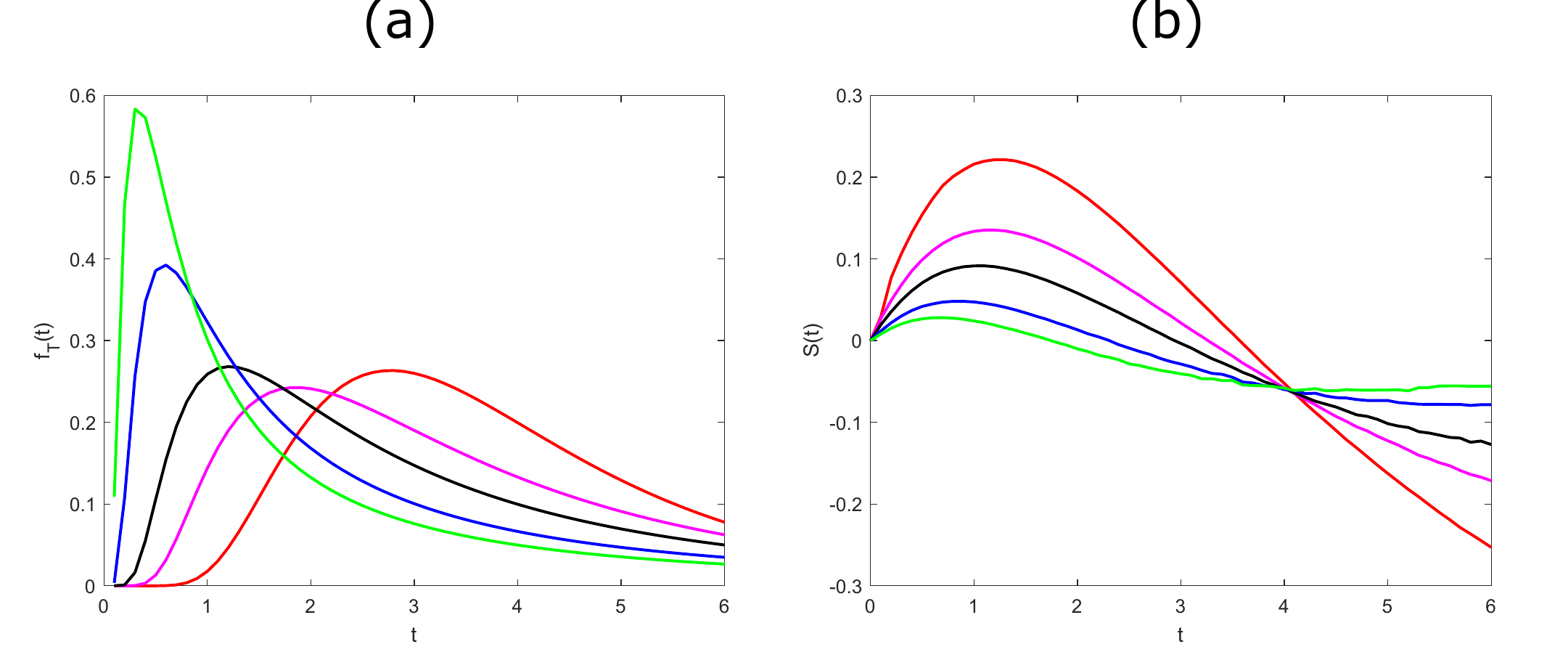}
	\caption{\label{FigIG}
	Probability densities (a) and evaluated boundaries (b) in the case of Inverse Gaussian-distributed FPTs with $\mathbb{E}(T)=4$. Different lines correspond to different shapes of the Inverse Gaussian densities,
$CV = 0.5$ (red), $CV = 0.75$ (magenta), $CV = 1$ (black), $CV = 1.5$
(blue), $CV = 2$ (green). The parameters of the two compartment model are $\alpha=0.33$,
$\beta=0.2$, $\mu=0$ and $\sigma = 1$.}
\end{figure}

Throughout all this paper, when not differently specified, we fix the values of the two compartment model as follows: $\alpha=0.33, \beta=0.2$ and $\sigma=1$.
We choose the shape of the IG distribution by fixing its mean  $\mathbb{E}[T]=4$ and different values of $CV$ (see Figure \ref{FigIG}, panel (a)). 
From (\ref{CV_IG}) we see that changes of $CV$ imply changes of the shape parameter $\lambda$. Moreover, as $CV$ increases, the density becomes more peaked. 
The corresponding shapes of the time varying thresholds are illustrated in panels (b) where the shapes of the boundary present a maximum that tends to disappear as $CV$ grows to higher values. 

\begin{figure}[t!]
\includegraphics[width=1\textwidth]{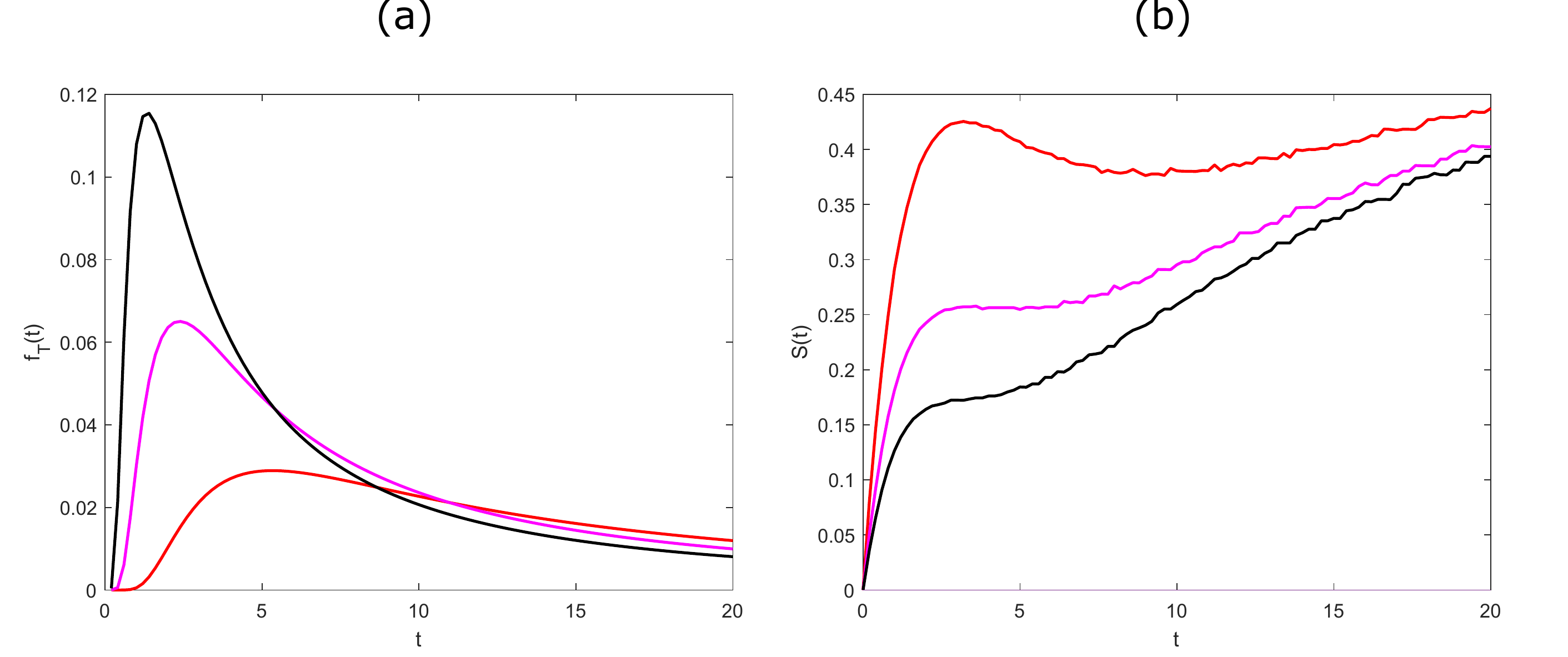}
	\caption{\label{FigIG W}
	Probability densities (a) and evaluated boundaries (b) in the case of Inverse Gaussian-distributed FPTs with heavy tails. Different lines correspond to different values of the parameter: $\lambda=16$ (red), $\lambda=7.11$ (magenta), $\lambda=4$ (black). The parameters of the two compartment model are $\alpha=0.33$, $\beta=0.2$, $\mu=0$ and $\sigma = 1$. Probability mass=[0.37 0.55 0.65] for $\lambda=[16 \text{ }7.11\text{ } 4]$.}
\end{figure}

When $\rho$ is finite, the IG distribution has light tails but if $\rho\rightarrow \infty$ the IG becomes 
\begin{equation}\label{IG Wiener}
f_T(t)=\left[ \frac{\lambda}{2\pi t^3}\right]^{1/2} \exp\left[ -\frac{\lambda}{2 t}\right], \quad \quad t\geq 0, 
\end{equation}
and it catches the heavy tails feature of interest for the analysis of some data \cite{GM,KYP,TIF}. The density (\ref{IG Wiener}) is known to be the density of the FPT of a Brownian motion with zero drift and diffusion coefficient $\nu$ through a constant boundary $b$ with the relation
\begin{equation}
\lambda=\frac{b^2}{\nu^2}.
\end{equation} 
Since $\mathbb{E}(T)=\infty$, it makes no sense to compute the $CV$ but, in order to compare light and heavy tails distributions, we use the same values of $\lambda$ in Figure \ref{FigIG} and \ref{FigIG W}. In Figure \ref{FigIG W}, different shapes of the pdfs (panel (a)) and of the corresponding boundaries (panel (b)) are shown. The heavy tails of this distribution determine a new shape for the threshold that has a decreasing maximum as $\lambda$ increases, followed by a minimum and by an increasing shape of the boundary. The maximum tends to disappear for large values of $\lambda$ and the values of the boundary are essentially positive. 
The growth of the boundary, for larger values of $t$, stop to allow the crossing of the samples determining the tail of the distribution. Figure \ref{FigIG W} refers to the time interval $[0,20]$, corresponding to a low probabilistic mass. A check for longer intervals does not change the results from a qualitative viewpoint, while higher probability masses are reached (figure not shown).



\subsection{Gamma random variable}\label{Sec:Gamma}
A random variable $T$ is Gamma distributed if its pdf is
\begin{equation}\label{f_T}
f_T(t)=\frac{\gamma ^\kappa }{\Gamma(\kappa)}t^{\kappa-1}e^{-\gamma t}, \quad \quad t\geq 0.
\end{equation}
Here, $\gamma>0$ is the rate parameter and $\kappa>0$ is the shape parameter.
Such a random variable is characterized by the following mean, variance and coefficient of variation 
\begin{eqnarray}
E(T)&=&\frac{\kappa}{\gamma}\\\nonumber
Var(T)&=&\frac{\kappa}{\gamma^2}\label{CV_Gamma}\\
CV&=&CV(T)=\frac{1}{\sqrt{\kappa}}.\nonumber
\end{eqnarray}

\begin{figure}[t!]
\includegraphics[width=1\textwidth]{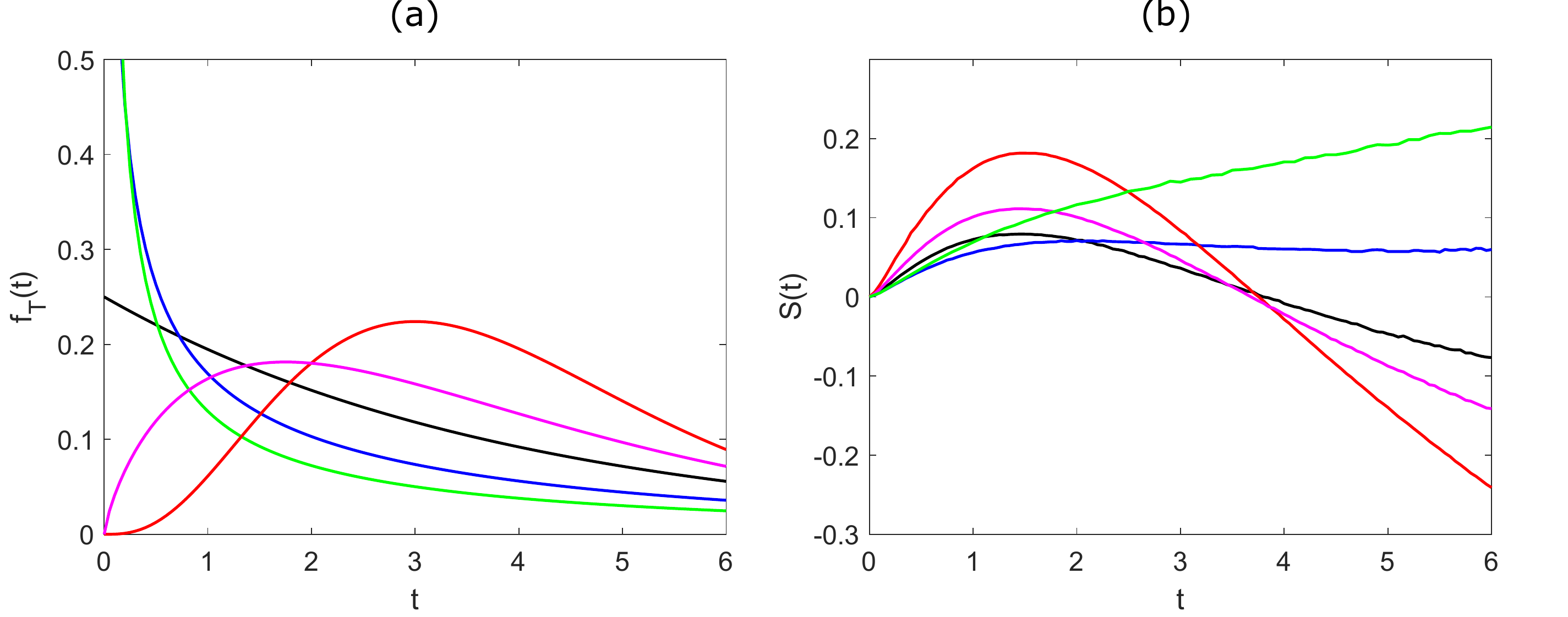}
	\caption{\label{Fig.G} Probability densities (a) and evaluated boundaries (b) in the case of Gamma-distributed FPTs with $\mathbb{E}(T)=4$. Different lines correspond to different shapes of the Gamma densities,
$CV = 0.5$ (red), $CV = 0.75$ (magenta), $CV = 1$ (black), $CV = 1.5$
(blue), $CV = 2$ (green). The parameters of the two compartment model are $\alpha=0.33$,
$\beta=0.2$, $\mu=0$ and $\sigma = 1$.}
\end{figure}

The shapes of Gamma pdf for different values of the parameters $\gamma$ and $\kappa$ can be seen in Figure \ref{Fig.G}, panel (a). 
The shapes of the Gamma pdf strongly change with the value of $CV$. We recall that $CV=1$ corresponds to the exponential distribution. Tails of the Gamma distribution are light, decaying to zero as an exponential. 
The corresponding shapes of the time varying thresholds are shown in panel (b) where the  values of the two compartment model are: $\alpha=0.33$, $\beta=0.2$ and $\sigma=1$. Here, as $CV$ increases, the maximum of the boundary disappears and the threshold time varying threshold becomes flat or, eventually when $CV=2$, increasing. This is the main difference of the boundary behavior with respect to the IG case.

\subsection{Comparison}
Often, IG and Gamma distributions appear as output of models of the same phenomenon but, for different choices of diffusion parameters. Hence, it seems useful to compare these distributions in terms of corresponding boundaries, using the Inverse FPT method.
Hence, in this subsection we compare the boundaries corresponding to IG and to Gamma distributions when mean and $CV$ are the same.  

\begin{figure}[t!]
\includegraphics[width=1\textwidth]{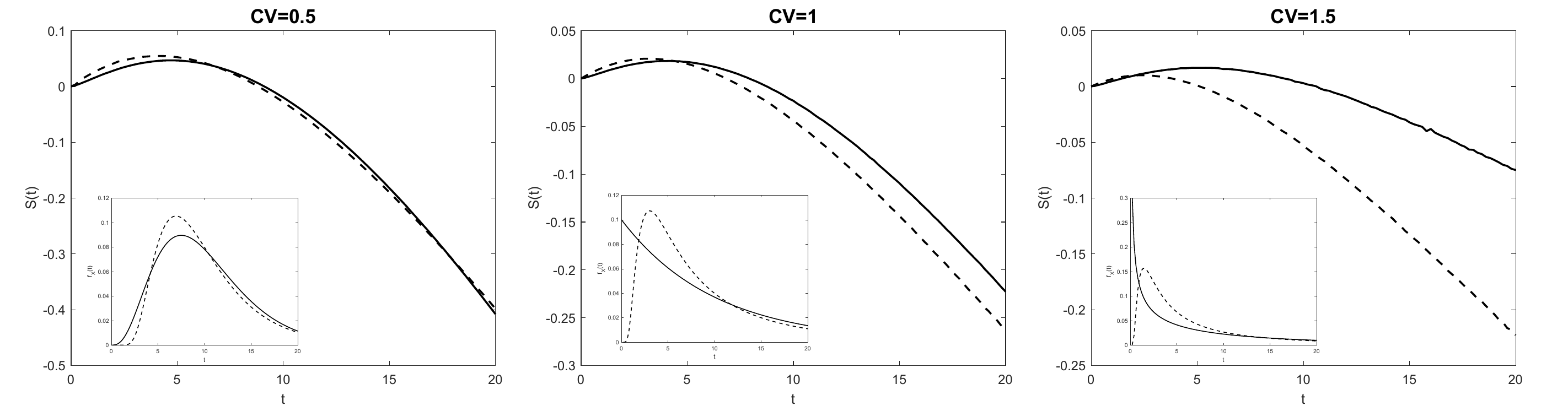}
	\caption{\label{FigGamGI} Comparison between the time varying boundaries corresponding to the Inverse Gaussian (dashed) and the Gamma-distributed (solid) FPTs with $\mathbb{E}(T)=10$ and with varying $CV$, the probability densities are on the sub-plots. The parameters of the two compartment model are $\alpha=0.02$,
$\beta=0.02$, $\mu=0$ and $\sigma = 0.4$.}
\end{figure}

Figure \ref{FigGamGI} shows the time varying boundaries corresponding to the IG (dashed) and the Gamma-distributed (solid) interspike intervals (ISIs) with the same mean value and the same $CV$. In this example $\mathbb{E}[T]=10$ while $CV=[0.5, 1, 1.5]$. Densities and corresponding boundaries become more and more different as we increase the $CV$ value (cf. inbox of Figure \ref{FigGamGI}). The different spreading of probability mass of the two classes of distributions is reflected in different shapes of the corresponding boundaries.  Since the IG density has heavier tails, the probability mass should not be consumed for short times. For this reason the boundary increases allowing crossings for large times.

\begin{figure}[t!]
\includegraphics[width=1\textwidth]{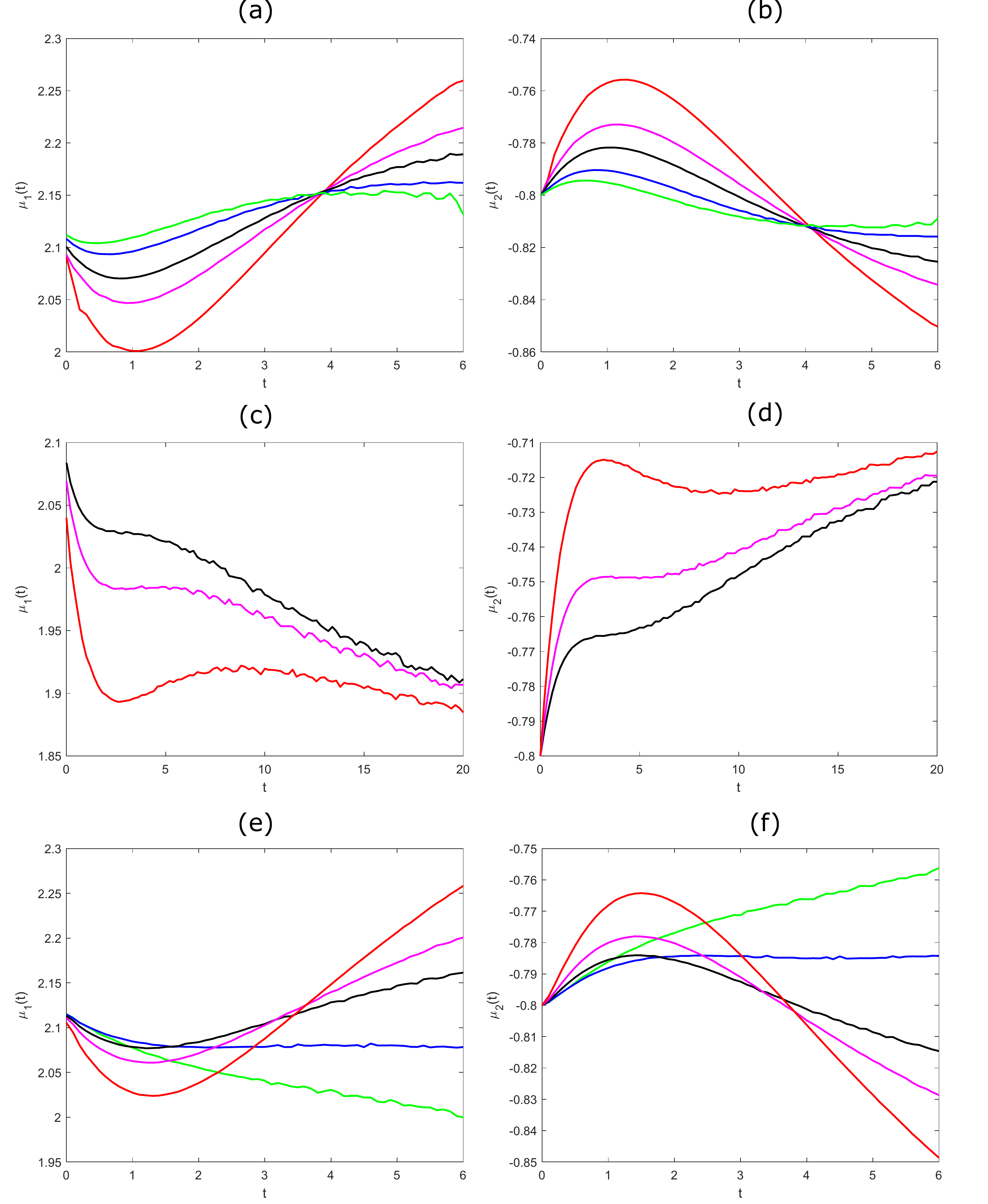}
	\caption{\label{FigG_IG_mu} Comparison between the values of the two components of $M(t)$ (\ref{M}) when the boundary is constant $\Sigma=4$ in the case of FPT distributed as IG (a-b), IG with heavy tails (c-d) and Gamma (e-f), respectively. The parameters and the colors of the functions are the same of the examples of Figures \ref{FigIG}, \ref{FigIG W} and \ref{Fig.G}. The parameters of the two compartment model are $\alpha=0.33$,
$\beta=0.2$, $\mu=0$ and $\sigma = 1$.}
\end{figure}

To help a physical interpretation of the results in terms of input of a two compartment model, in Figure \ref{FigG_IG_mu}  we compare the behavior of the two components (\ref{M}) of $M(t)$, when the boundary is transformed into a constant through (\ref{trasformazione}). We ask what would be the input to both compartments if the output distribution is fixed and threshold is constant. We illustrate the cases of FPT distributed as IG (a-b), IG with heavy tails (c-d) and Gamma (e-f), respectively. Reinterpreting the time dependent boundary in terms of a modification of the drift allows to interpret our results in terms of increasing or decreasing drift. We note that a positive drift on the first component is always necessary to obtain the prescribed FPT distribution. When $CV\leq 1$,  FPTs distributed as Gamma or IG imply similar input. On the contrary, when tails of IG are heavy (panel (c)) or $CV$ is large enough (panels (a) or (e)), the drift term of the first component strongly changes becoming decreasing. Interestingly the behavior of the second component $\mu_2(t)$ (panels (b),(d) and (f)) is the opposite of that of $\mu_1(t)$. 

\begin{figure}[t!]
\includegraphics[width=1\textwidth]{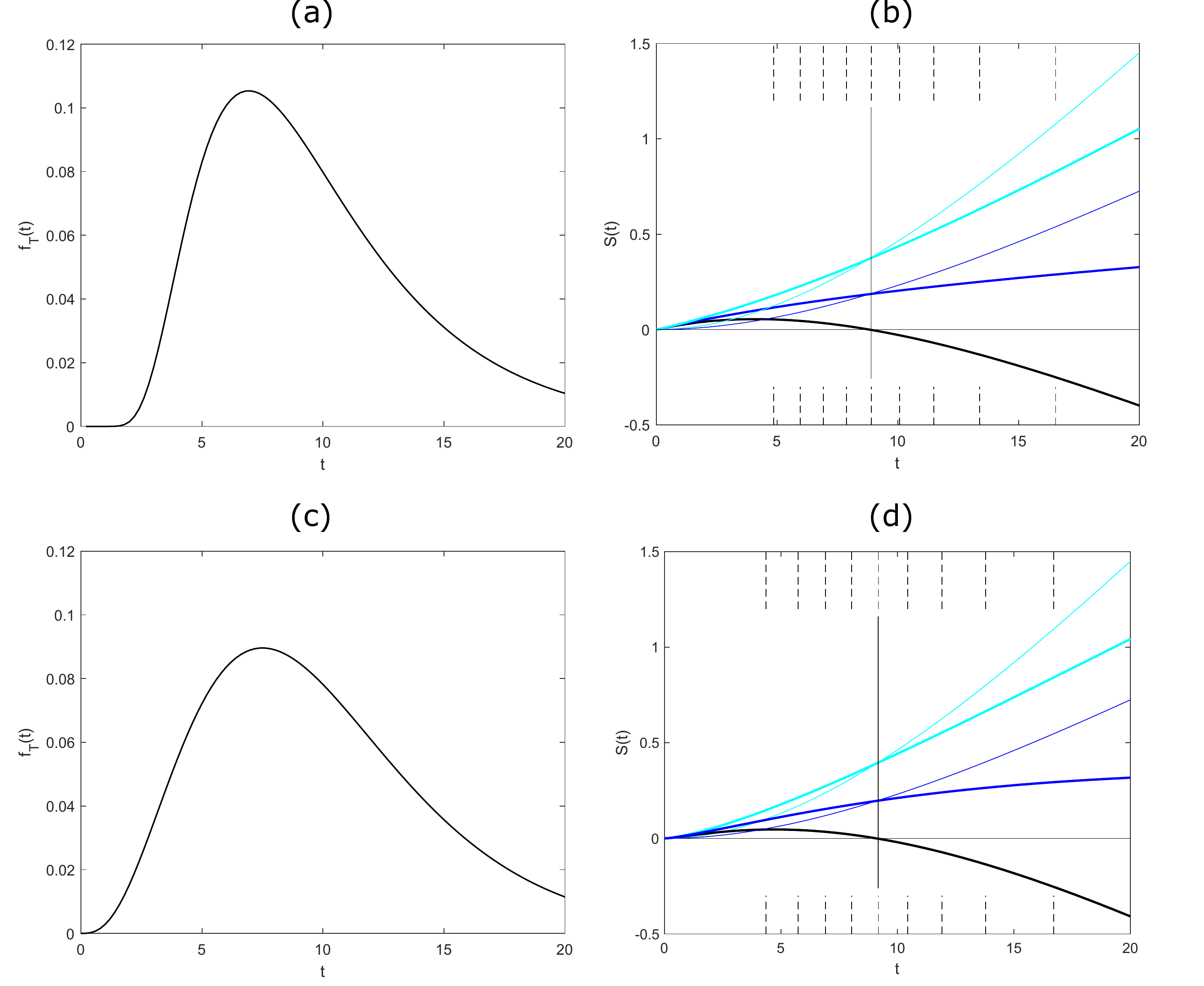}
	\caption{\label{Fig mu CV=0.5}Inverse-Gaussian (a) and Gamma (c) distributed FPTs with mean FPT equal to 10 and $CV=0.5$. Corresponding boundaries and mean value of the first component (b,d). Different lines correspond to different values of the mean input: $\mu=0$ (black), $\mu = 0.3$ (blue), $\mu = 0.6$ (cyan). Thicker lines correspond to the time varying boundaries.
The vertical dotted lines give the FPT distribution quantiles. The parameters of the two compartment model are  $\alpha=0.02$,
$\beta=0.02$, $\sigma = 0.4$, while $\mu$ varies as specified above.}
\end{figure}

\begin{figure}[t!]
\includegraphics[width=1\textwidth]{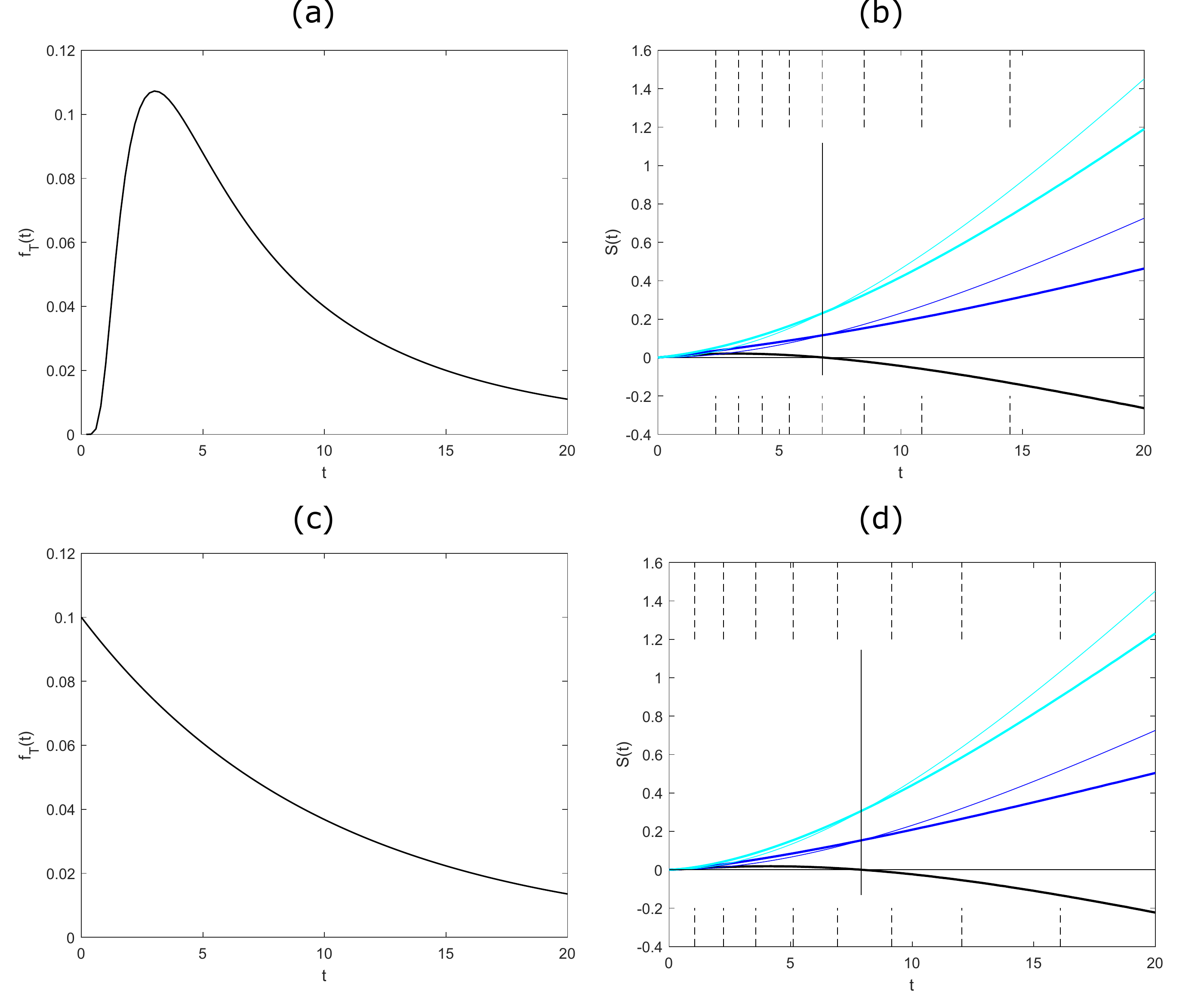}
\caption{\label{Fig mu CV=1}Inverse-Gaussian (a) and Gamma (c) distributed FPTs with mean FPT equal to 10 and $CV=1$. Corresponding boundaries and mean value of the first component (b,d). Different lines correspond to different values of the mean input: $\mu=0$ (black), $\mu = 0.3$ (blue), $\mu = 0.6$ (cyan). Thicker lines correspond to the time varying boundaries.
The vertical dotted lines give the FPT distribution quantiles. The parameters of the two compartment model are  $\alpha=0.02$,
$\beta=0.02$, $\sigma = 0.4$, while $\mu$ varies as specified above.}
\end{figure}

Lastly, we change the comparison criterion and we apply the Inverse FPT method varying the values of the parameter $\mu$ in (\ref{model}). We fix the values of the model as follows: $\alpha=0.02$, $\beta=0.02$ and  $\sigma=0.4$. We consider examples of boundaries corresponding IG or Gamma spiking densities for $CV=0.5$ (Figure \ref{Fig mu CV=0.5}) or $CV=1$ (Figure \ref{Fig mu CV=1}).
 
In the figures, we also compare the boundaries (thicker line) with the mean of the first component $\mathbb{E}[X_1(t)]$. 
We note that boundary always intersects the function $\mathbb{E}[X_1(t)]$. Interestingly, if we fix the firing FPT and its $CV$, the intersection value is the same for different values of the parameter $\mu$. This fact can be easily understood by noting that a change in $\mu$ determines the same shift both on the mean value of $X_1(t)$ and on the boundary. However, this value changes considering different $CV$s.

\section{Application to neuroscience}\label{Sec: Results}

We give here an example of application of the Inverse FPT method to neuroscience. 

Simplest neuronal models resort to one-dimensional processes to describe the membrane potential evolution. This choice implies a strong simplification of the neuronal structure that is identified by a single point.
 More complex models introduce bivariate stochastic processes to discriminate the membrane potential dynamics in the dendritic or in the trigger zone \cite{LR}. Neurophysiological reasons suggest the existence of an interaction between the membrane potential dynamics in the two zones and, when the first component (the trigger one) attains a boundary value, the neuron releases a spike. A reasonable simplification allows to add a noisy term only to the dendritic component.  
The  reset after the spike can include both the components or only the trigger zone. Here, we will consider only the case of total resetting of both components to a resting value that we fix, for simplicity, equal to zero.\\

\begin{figure}[!h]
		\includegraphics[width=1\textwidth]{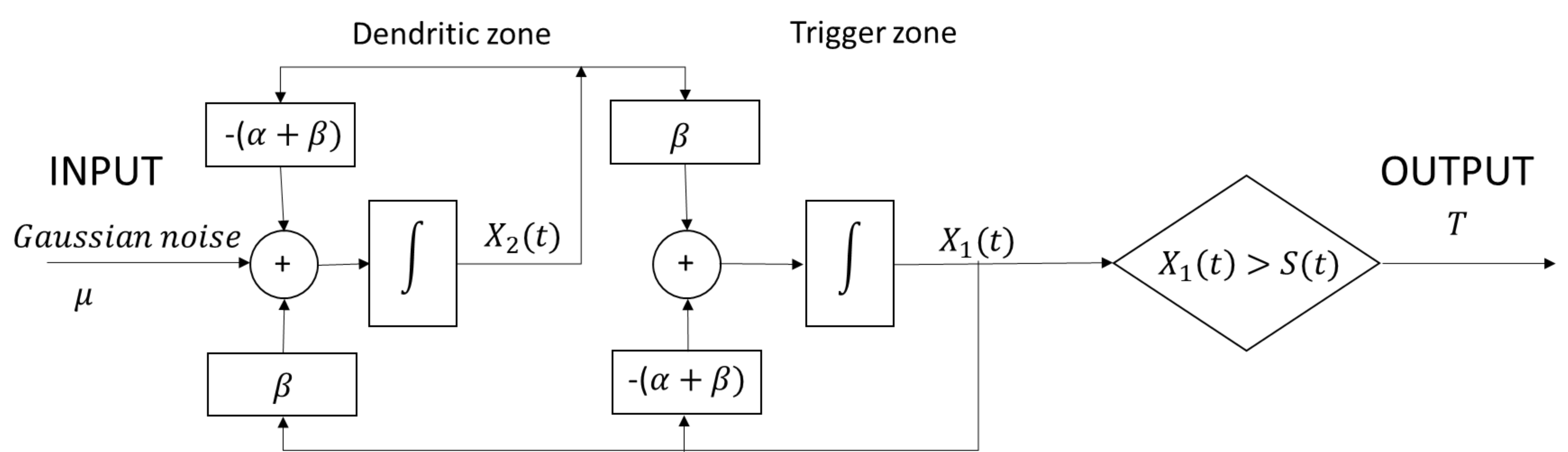}
	\caption{Schematic representation of the two-compartment model}
	\label{BiModel}
\end{figure}

In this framework, the stochastic process (\ref{model}) describes the depolarization of the trigger zone and the dendritic one, respectively \cite{LR}. The model assumes that external inputs, with intensity $\mu$ and variability $\sigma$, influence the second compartment and a weight $\beta$ takes into account the interconnection between the parts of the neuron. Moreover, the constant $\alpha>0$ accounts for the spontaneous membrane potential decay (cf. Figure \ref{BiModel}).
Then, the FPT $T$ mimics the ISI of the neuron and the boundary $S(t)$ corresponds to the spiking threshold for the neuron.

Often, IG and Gamma distribution fit neuronal data and FPT may help to interpret the presence of these distributions. Here, we reinterpret Figures \ref{FigIG}-\ref{Fig mu CV=1} in the neuronal model framework. Hence, constants  
 $\alpha$ and $\beta$ will be measured in $\text{ ms}^{-1}$, while $\mu$ will be measured in $\text{mV ms}^{-1}$ and $\sigma$ in $\text{ mV ms}^{-1/2}$. 
 
Figures \ref{FigIG}, \ref{FigIG W}, \ref{Fig.G} and \ref{FigGamGI} reinterpret the heaviness of the tail of the ISI distributions in terms of the threshold shapes. As we increase the $CV$ value, IG and Gamma densities and the corresponding boundaries become more and more different. This means that $CV$ plays an important role in the formulation of the model. Moreover, in the case of the Gamma ISI distribution, the slope of the threshold becomes increasing when CV is large enough. A similar increasing behavior of the boundary could be obtained with IG distributed ISIs with heavy tails (cf. Figure \ref{FigIG W}).

In Figures \ref{Fig mu CV=0.5} and \ref{Fig mu CV=1} we investigate not only the behavior of the time varying firing threshold but also the dynamics of the underlying two compartment neuronal model. 
The mean membrane potentials and the corresponding boundaries are plotted for different values of the mean input $\mu$ and for different values of $CV$.
For low input $\mu$, the curves exhibit a maximum after which the firing threshold starts to decrease. As the input $\mu$ increases, the firing threshold from concave and decreasing become convex and increasing. Indeed, a big input $\mu$ facilitates the spiking. Therefore, to obtain the assigned distribution, the threshold must move away, becoming increasing.

\begin{figure}[t!]
\includegraphics[width=1\textwidth]{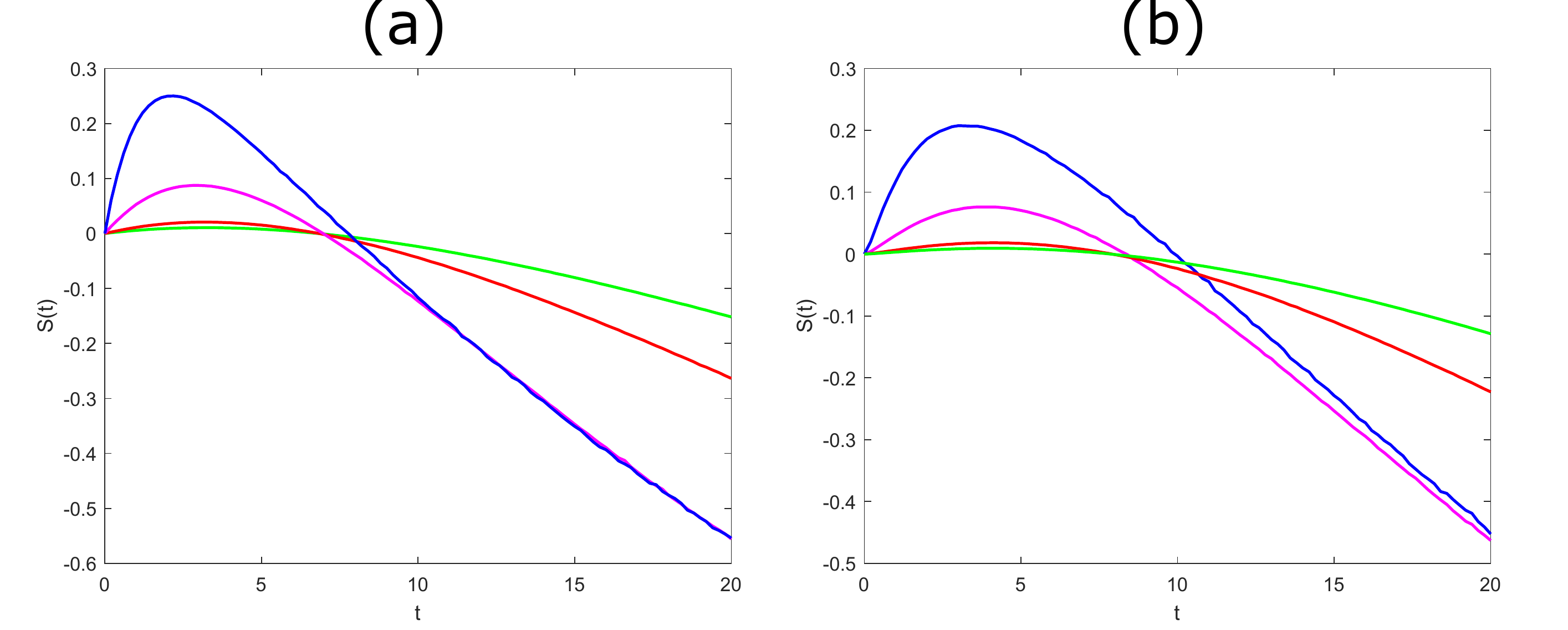}
	\caption{\label{Fig. Comparison Beta}
	Boundaries corresponding to Inverse Gaussian (a) and Gamma (b) distributed FPTs with $\mathbb{E}(T)=10$ and $CV=1$. Different lines correspond to different values of $\beta$: $\beta=0.01$ (green), $\beta = 0.02$ (red), $\beta=0.1$ (magenta), $\beta=0.5$ (blue).The parameters of the two compartment model are  $\alpha=0.02$,
$\mu=0$, $\sigma = 0.4$, while $\beta$ varies as specified above.}
\end{figure}

Lastly in Figure \ref{Fig. Comparison Beta} we study the role of the parameter $\beta$ in the model, applying the Inverse FPT method to IG (panel (a)) and Gamma (panel (b)) FPT distribution and varying the values of the parameter $\beta$. As $\beta$ decreases, the boundary becomes almost constant and equal to zero. This is consistent with the fact that, for $\beta = 0$, the two components of the process gets independent and $X_1$ is deterministic and equal to zero, since $X(0)=0$. Then, in order to have a crossing and to get a prescribed distribution, the threshold should approach zero.

\section{Conclusions}
The extension of the Inverse FPT method to two-dimensional OU diffusion processes allows to study the shape of the boundaries for a  given FPT pdf. 
We applied the algorithm to FPT distributed as an Inverse Gaussian and as a Gamma random variable. Differences in the boundary shape corresponding to FPTs with heavy or light tails enlighten different features of the corresponding two compartment model.  

Lastly, we reinterpret the obtained results in a neuroscience framework. The shape of the boundaries corresponding to different firing distributions may enlighten features of the model eventually recognizing instances of scarce physiological significance such as diverging thresholds.

\section*{Acknowledgement}
The authors are grateful to Professor L. Sacerdote and Professor L. Kostal for their interesting and useful comments and suggestions.

\end{document}